\documentclass{article}
\usepackage{latexsym}
\usepackage{amsfonts}
\newtheorem{theorem}{Theorem}
\newtheorem{lemma}{Lemma}
\newtheorem{proposition}{Proposition}

\begin{document}
\def\firstone{``first $1$''}
\def\firstnegone{``first $-1$''}
\def\tb{{\tilde b}}
\def\grpz{{\mathbb Z}}
\def\qed{{\hfill $\Box$}}
\title{A lower bound for the Chung-Diaconis-Graham random process}
\author{Martin Hildebrand\\Department of Mathematics and Statistics\\
University at Albany\\State University of New York\\Albany, NY 12222}

\maketitle

\begin{abstract}
Chung, Diaconis, and Graham considered random processes of the form
$X_{n+1}=a_nX_n+b_n\pmod p$ where $p$ is odd, $X_0=0$, $a_n=2$ always,
and $b_n$ are
i.i.d. for $n=0,1,2,\dots$. In this paper, we show that if
$P(b_n=-1)=P(b_n=0)=P(b_n=1)=1/3$, then there exists a constant $c>1$
such that $c\log_2p$ steps are not enough to make $X_n$ get close
to uniformly distributed on the integers mod $p$.

\end{abstract}

\section{Introduction}
In \cite{cdg}, Chung, Diaconis, and Graham considered random processes
of the form
\[
X_{n+1}=2X_n+b_n\pmod p
\]
where $X_0=0$, $b_n$ are i.i.d. for $n=0,1,2,\dots$, and $p$ is odd.
They focussed on the case where $P(b_n=1)=P(b_n=0)=P(b_n=-1)=1/3$.
These random processes have some similarity to certain pseudorandom
sequences used by computers. Subsequently some generalizations of this random
process have been considered. See, for example, \cite{asci}, \cite{mvhannprob},
\cite{mvhima}, and \cite{mvhecp}. Suppose $P_n(s)=P(X_n=s)$ where $s\in
\grpz/p\grpz$. Define the variation distance of a probability
$P$ on $\grpz/p\grpz$ from the uniform distribution $U$ on $\grpz/p\grpz$
by
\[
\|P-U\|={1\over 2}\sum_{s\in\grpz/p\grpz}|P(s)-1/p|=
\max_{A\subseteq\grpz/p\grpz}|P(A)-U(A)|.
\]
They showed that for almost all odd $p$, if
$N\ge{\log p\over \log(9/5)}+c$, then $\|P_N-U\|=O((5/9)^c)$.
They also state that more complicated arguments give the following
result: For any $\epsilon>0$ and almost all odd $p$, if
$N\ge(\hat c+\epsilon)\log_2p$, then $\|P_N-U\|<\epsilon$ where
\[
\hat c=\left(1-\log_2\left({5+\sqrt{17} \over 9}\right)\right)^{-1}=
1.01999186\dots
\]

Note that the values of $X_N$, when viewed as integers, range between
$-2^N+1$ and $2^N-1$, inclusive. Thus if $N<(1-\epsilon)\log_2p$ where
$\epsilon>0$ is given, the
number of values in this range is at most $2p^{1-\epsilon}-1$ and $\|P_N-U\|
>1-(2p^{1-\epsilon}-1)/p\rightarrow 1$ as $p\rightarrow\infty$.

Chung, Diaconis, and Graham~\cite{cdg} speculate, ``It is conceivable that
in fact $(1+o(1))\log_2p$ steps are enough for almost all [odd] $p$
to force $P_N$ to converge to uniform.'' However, we shall show that
this statement, although described as conceivable in \cite{cdg}, in fact is
false. In particular, we shall show the following theorem:
\begin{theorem}
\label{mainthm}
If $P(b_n=1)=P(b_n=0)=P(b_n=-1)=1/3$ and $X_n$ and $P_n$ are as above,
then there exists a value $c_1>1$ such that if $n=n(p)<c_1\log_2p$, then
$\|P_n-U\|\rightarrow 1$ as $p\rightarrow\infty$.
\end{theorem}

To motivate somewhat the proof of this theorem, we shall also prove
the following theorem:
\begin{theorem}
\label{warmup}
If $P(b_n=1)=0.4$, $P(b_n=0)=0.6$, and $X_n$ and $P_n$ are as above, then 
there exists a value $c_2>1$ such that if $n=n(p)<c_2\log_2p$, then
$\|P_n-U\|\rightarrow 1$ as $p\rightarrow\infty$.
\end{theorem}

\section{Proof of Theorem~\protect\ref{warmup}}

First observe the following proposition:
\begin{proposition}
\label{propone}
If $X_0=0$ and $X_{n+1}=2X_n+b_n$ for $n\ge 0$, then
\[
X_n=\sum_{i=0}^{n-1}2^{n-1-i}b_i.
\]
\end{proposition}

Now suppose $P(b_n=1)=0.4$ and $P(b_n=0)=0.6$. Let $A_n=|\{m:0\le m\le n-1,
b_m=1\}|$. By elementary arguments, for any $\epsilon>0$, $P((0.4-\epsilon)n
<A_n<(0.4+\epsilon)n)\rightarrow 1$ as $n\rightarrow\infty$. 
Thus, except on a set which has
probability approaching $0$ as $n\rightarrow\infty$, $X_n$ takes on at
most
\[
\sum_{j=\lceil(0.4-\epsilon)n\rceil}^{\lfloor (0.4+\epsilon)n\rfloor}
{n\choose j}
\]
different values. We shall assume that $0.4+\epsilon<0.5$. Note that
Stirling's formula implies
\begin{eqnarray*}
&&\sum_{j=\lceil(0.4-\epsilon)n\rceil}^{\lfloor (0.4+\epsilon)n\rfloor}
{n\choose j}\\
&\le&
(2\epsilon n+1){n\choose \lfloor(0.4+\epsilon)n\rfloor}\\
&\le&{(2\epsilon n+1)c_3
n^n\sqrt{2\pi n} \over
((0.4+\epsilon)n)^{(0.4+\epsilon)n}
((0.6-\epsilon)n)^{(0.6-\epsilon)n}
2\pi\sqrt{((0.4+\epsilon)n)((0.6-\epsilon)n)}}\\
&\le&
{c_3(2\epsilon n+1)\sqrt{2\pi n}\over 2\pi\sqrt{((0.4+\epsilon)n)
((0.6-\epsilon)n)}}\cdot
{1\over (0.4+\epsilon)^{(0.4+\epsilon)n}(0.6-\epsilon)^{(0.6-\epsilon)n}}
\end{eqnarray*}
where $c_3$ is a positive constant.
Note that
\[
(0.4+\epsilon)^{(0.4+\epsilon)n}(0.6-\epsilon)^{(0.6-\epsilon)n}
=
2^{n((0.4+\epsilon)\log_2(0.4+\epsilon)+(0.6-\epsilon)\log_2(0.6-\epsilon))}.
\]
It can be shown that if $0<\epsilon<0.1$, then
\[
-((0.4+\epsilon)\log_2(0.4+\epsilon)+(0.6-\epsilon)\log_2(0.6-\epsilon))<1.
\]
If
\[
c_2<{1\over
-((0.4+\epsilon)\log_2(0.4+\epsilon)+(0.6-\epsilon)\log_2(0.6-\epsilon))}
\]
and $n=n(p)<c_2\log_2p$, then
\[
{\sum_{j=\lceil(0.4-\epsilon)n\rceil}^{\lfloor (0.4+\epsilon)n\rfloor}
{n\choose j}\over p}\rightarrow 0
\]
as $p\rightarrow\infty$. Thus $\|P_n-U\|\rightarrow 1$ as $p\rightarrow\infty$
if $n=n(p)<c_2\log_2p$. Note that $c_2$ can be chosen so that $c_2>1$ since
\[{1\over
-((0.4+\epsilon)\log_2(0.4+\epsilon)+(0.6-\epsilon)\log_2(0.6-\epsilon))}>1.\]
\qed

\section{Overview of the Proof of Theorem~\protect\ref{mainthm}}

By Proposition~\ref{propone}, $X_n$ is determined by the $n$-tuple
$(b_0, b_1, \dots, b_{n-1})$. However, if $P(b_n=1)=P(b_n=0)=P(b_n=-1)$,
then many possible $n$-tuples $(b_0, b_1, \dots, b_{n-1})$ may give the
same value for $X_n$. For example, if $n=3$, the $3$-tuples $(1,-1,-1)$,
$(0,1,-1)$, and $(0,0,1)$ all give $X_3=1$. We shall place this $n$-tuple
in a standard form $(\tb_0,\tb_1,\dots,\tb_{n-1})$ so that
\[
2^{n-1}\tb_0+2^{n-2}\tb_1+\dots+\tb_{n-1}=
2^{n-1}b_0+2^{n-2}b_1+\dots+b_{n-1}.
\]
In this standard form, either none of $\tb_0$, $\tb_1$, $\dots$, $\tb_{n-1}$
are $-1$ or none of $\tb_0$, $\tb_1$, $\dots$, $\tb_{n-1}$ are $1$.
In the first case (excluding the event where $b_0, b_1, \dots, b_{n-1}$ are
all $0$), we shall show that for every $\epsilon>0$, the number of values
$a$ in $\{1,\dots,n-1\}$ such that both $\tb_{a-1}$ and $\tb_a$ are 
$1$ lies between $(4/18-\epsilon)n$ and $(4/18+\epsilon)n$ except for
events which have probability approaching $0$ as $n\rightarrow\infty$.
Likewise, in the second case (excluding the event where $b_0, b_1, \dots,
b_{n-1}$ are all $0$), the number of values $a$ in $\{1,2,\dots,n-1\}$ such
that both $\tb_{a-1}$ and $\tb_a$ are both $-1$ lies between 
$((4/18)-\epsilon)n$
and $((4/18)+\epsilon)n$ except for events which have probability approaching
$0$ as $n\rightarrow\infty$. A Stirling's formula argument will give the
theorem.

We shall divide the $n$-tuples $(b_0, b_1, \dots, b_{n-1})$ into three cases.
In the first case, there exists a value $j$ in $\{0, 1, \dots, n-1\}$ such that
$b_j=1$ and $b_k=0$ if $0\le k<j$. We call this case \firstone. In the
second case, there exists a value $j$ in $\{0, 1, \dots, n-1\}$ such that
$b_j=-1$ and $b_k=0$ if $0\le k<j$. We call this case \firstnegone.
In the third case, $b_0=b_1=\dots=b_{n-1}=0$. As $n\rightarrow\infty$,
the probability of \firstone approaches $1/2$, and the probability of
\firstnegone approaches $1/2$. (Both probabilities are $(1/2)(1-(1/3)^n)$.)
In \firstone, none of $\tb_0$, $\tb_1$, $\dots$, $\tb_{n-1}$ are $-1$. In
\firstnegone, none of $\tb_0$, $\tb_1$, $\dots$, $\tb_{n-1}$ are $1$.
We shall give detailed arguments for \firstone; the arguments for \firstnegone
are similar.

If we are in \firstone, consider the infinite sequence $(b_0,b_1,b_2,\dots)$.
After some leading zeroes, with probability $1$, this sequence consists of
strings of ``blocks'' $B_1, B_2, \dots$ where $B_i$ is a finite string 
starting with $1$ and having no other $1$'s in it. 
For example, if $(b_0,b_1,\dots,b_{10})=(0,0,1,-1,0,1,0,1,-1,1,1)$, then
there are two leading zeroes, $B_1=(1,-1,0)$, $B_2=(1,0)$, $B_3=(1,-1)$,
and $B_4=(1)$.
Also we say that the first coordinate of $B_1$ is $b_2$ and that $B_1$ has
$b_2$, $b_3$, and $b_4$ as its coordinates. Note that the blocks
$B_1, B_2, \dots$ are i.i.d. given that we are in \firstone. 

Given an infinite series $(b_0,b_1,b_2,\dots)$, technically
the values $\tb_0,\tb_1,\dots,\tb_{n-1}$ may depend on
$n$. For example, if $(b_0,b_1,b_2)=(0,0,1)$, then $(\tb_0,\tb_1,\tb_2)=
(0,0,1)$ and $\tb_2=1$ if $n=3$. If 
$(b_0,b_1,\dots,b_{10})=(0,0,1,-1,0,1,0,1,-1,1,1)$, then
$(\tb_0,\tb_1,\dots,\tb_{10})=(0,0,0,1,0,1,0,0,1,1,1)$ and $\tb_2=0$ if $n=11$
even though $b_0$, $b_1$, and $b_2$ are unchanged. Suppose $a$ is
such that $B_i$ has $b_a$ as one of its coordinates. Then for all $n>a$ such
that $B_i$ does not also have $b_n$ as one its coordinates,
$(\tb_0,\tb_1,\dots,\tb_a)$ will no longer vary with $n$.

\section{Number of $a$ such that $\tb_{a-1}=\tb_a=1$}

We shall consider several distinct ways to get values $a$ such that
$\tb_{a-1}=\tb_a=1$. These ways are detailed in the following lemmas.

\begin{lemma}
\label{lemmaone}
Let $n_1$ be the number of $a$ in $\{1,\dots,n-1\}$ such that
$\tb_{a-1}=1$, $\tb_a=1$, $b_{a-1}=1$, and $b_a=1$. Let $\epsilon>0$
be given. Given that we are in \firstone, the probability that
$((1/18)-\epsilon)n<n_1<((1/18)+\epsilon)n$ approaches $1$ as
$n\rightarrow\infty$.
\end{lemma}

\begin{lemma}
\label{lemmatwo}
Let $n_2$ be the number of $a$ in $\{1,\dots,n-1\}$ such that
$\tb_{a-1}=1$, $\tb_a=1$, $b_{a-1}=1$, and $b_a\ne 1$. Then $n_2=0$.
\end{lemma}

\begin{lemma}
\label{lemmathree}
Let $n_3$ be the number of $a$ in $\{1,\dots,n-1\}$ such that
$\tb_{a-1}=1$, $\tb_a=1$, $b_{a-1}\ne 1$, and $b_a=1$. Let $\epsilon>0$
be given. Given that we are in \firstone, the probability that
$((1/18)-\epsilon)n<n_3<((1/18)+\epsilon)n$ approaches $1$ as
$n\rightarrow\infty$.
\end{lemma}

\begin{lemma}
\label{lemmafour}
Let $n_4$ be the number of $a$ in $\{1,\dots,n-1\}$ such that
$\tb_{a-1}=1$, $\tb_a=1$, $b_{a-1}\ne 1$, and $b_a\ne 1$. Let $\epsilon>0$
be given. Given that we are in \firstone, the probability that
$((1/9)-\epsilon)n<n_4<((1/9)+\epsilon)n$ approaches $1$ as
$n\rightarrow\infty$.
\end{lemma}

{\it Proof of Lemma~\ref{lemmatwo}}: If we are not in \firstone, then 
$\tb_a$ is never $1$. If we are in \firstone, then one obtains
$(\tb_0, \tb_1, \dots, \tb_{n-1})$ from $(b_0, b_1, \dots, b_{n-1})$
as follows. If $j$ is such that $b_j=1$ and $b_k=0$ whenever $0\le k<j$, then
$\tb_k=0$ whenever $0\le k<j$. Otherwise, for each $j_0$ such that $b_{j_0}=1$,
let $j_1=\min(n,\min\{\ell:\ell>j_0, b_{\ell}=1\})$. (By convention, assume
that the minimum of an empty set is $\infty$.) If $b_k=0$ for all $k$ with
$j_0<k<j_1$, then $\tb_{j_0}=1$ and $\tb_k=0$ for all $k$ with $j_0<k<j_1$.
Otherwise $\tb_{j_0}=0$, and one can figure out the unique values for
$\tb_k$ in $\{0,1\}$ when $j_0<k<j_1$. Lemma~\ref{lemmatwo} follows.\qed

{\it Proof of Lemma~\ref{lemmafour}}: To prove Lemma~\ref{lemmafour}, let
$j$ be such that $b_j=1$ and $b_k=0$ whenever $0\le k<j$. Suppose that 
$a-1>j$ and $a<n$. Then $P(b_{a-1}\ne 1, b_a\ne 1)=4/9$. Suppose
$j_0<a<j_1$ with $b_{j_0}=1$ and 
$j_1=\min(n,\min\{\ell:\ell>j_0, b_{\ell}=1\})$.
Given $j_0$ and $j_1$, there are $2^{j_1-j_0-1}$ possibilities for
$(b_{j_0+1}, \dots, b_{j_1-1})$ and $2^{j_1-j_0-1}$ possibilities for
$(\tb_{j_0+1}, \dots, \tb_{j_1-1})$. The possibilities for
$(b_{j_0+1},\dots,b_{j_1-1})$, which range from $(0,\dots,0)$ to
$(-1,\dots,-1)$, are in one-to-one correspondence with the
possibilities for $(\tb_{j_0+1},\dots,\tb_{j_1-1})$, which range
from $(0,\dots,0)$ to $(1,\dots,1)$.
Thus $P(\tb_{a-1}=1, \tb_a=1 | b_{a-1}\ne 1, b_a\ne 1)=1/4$, and
$P(\tb_{a-1}=1, \tb_a=1, b_{a-1}\ne 1, b_a\ne 1)=1/9$. Let
$C_a=\{\tb_{a-1}=1, \tb_a=1, b_{a-1}\ne 1, b_a\ne 1\}$. Conditioned on $j$
such that $b_j=1$ and $b_k=0$ for $0\le k<j$, the events $C_a$ for
$a-1>j$, $a<n$, and $a$ even are independent, and the events $C_a$ for
$a-1>j$, $a<n$, and $a$ odd are independent. Since $P(j>\epsilon_1 n)\
\rightarrow 0$ as $n\rightarrow \infty$ (given that we are
in \firstone) for each $\epsilon_1>0$, Lemma~\ref{lemmafour}
follows by elementary arguments. \qed

{\it Proof of Lemma~\ref{lemmaone}}: To prove this lemma, suppose that we
are in \firstone and $b_{a-1}=1$ with $a$ in $\{1,\dots,n-1\}$. Then
$\tb_{a-1}=1$ and $\tb_a=1$ if and only if $b_a=1$ and $b_k=0$ for
all $k$ with $a<k<j_1$ where $j_1=\min(n,\min\{\ell>a:b_{\ell}=1\})$.
Let us consider the infinite sequence $(b_0, b_1, \dots)$. 
For positive integers $i$, let $D_i$ be the event that 
$B_i=(1)$ and $B_{i+1}$ has no $-1$'s in it.
Note that $P(D_i)=(1/3)(\sum_{i=1}^{\infty}(1/3)^i)=1/6$. Observe that
$D_1$, $D_3$, $D_5$, etc. are independent and that $D_2$, $D_4$, $D_6$, etc.
are independent. Furthermore, given $\epsilon_1>0$, with probability
approaching $1$ as $n\rightarrow\infty$, the number of $a$ in
$\{1,\dots,n-1\}$ such that $b_{a-1}=1$ lies between $((1/3)-\epsilon_1)n$
and $((1/3)+\epsilon_1)n$ given that we are in \firstone. 
Suppose we are given $\epsilon^{\prime}>0$. Choose $\epsilon_1>0$
so that $\epsilon_1<6\epsilon^{\prime}$. Then with probability
approaching $1$ as $n\rightarrow\infty$, at least $((1/18)-\epsilon^{\prime})n$
events $D_i$ occur with $i<((1/3)-\epsilon_1)n$ while at most
$((1/18)+\epsilon^{\prime})n$ events $D_i$ occur with 
$i\le((1/3)+\epsilon_1)n$.
Thus given
$\epsilon^{\prime}>0$, the number of $i$ such that $D_i$ occurs and
the first coordinate of $B_{i+1}$ is $b_{\ell}$ for some $\ell<n$
is, with probability approaching $1$ as $n\rightarrow\infty$,
between $((1/18)-\epsilon^{\prime})n$ and $((1/18)+\epsilon^{\prime})n$.
This number of $i$ is within $1$ of the number of $a$ in Lemma~\ref{lemmaone};
the only possible difference occurs 
when the block $B_{i+1}$ has $b_n$ as one of its coordinates.
\qed

{\it Proof of Lemma~\ref{lemmathree}}: To prove this lemma, suppose we
are in \firstone, $b_{a-1}\ne 1$, and $b_a=1$ with $a$ in
$\{1,\dots,n-1\}$. Then $\tb_{a-1}=1$ and $\tb_a=1$ if and only if
$b_{a-1}=-1$ and $b_k=0$ for all $k$ with
$a<k<j_1$ where $j_1=\min(n,\min\{\ell>a: b_{\ell}=1\})$. Let us
consider the infinite sequence $(b_0,b_1,\dots)$. 
For positive integers $i$, let
$E_i$ be the event that $B_{i+1}$ has no $-1$'s in it and that $B_i$ ends
with $-1$.
Note that $P(E_i)=1/6$, that
$E_1, E_3, E_5, \dots$ are independent, and that
$E_2, E_4, E_6, \dots$ are independent.
Furthermore, given $\epsilon_1>0$, with
probability approaching $1$ as $n\rightarrow\infty$, the number of $a$
in $\{1,\dots,n-1\}$ with $b_a=1$ lies between $((1/3)-\epsilon_1)n$
and $((1/3)+\epsilon_1)n$ given that we are in \firstone.
Thus given $\epsilon^{\prime}>0$, the number of $i$ such that $E_i$ occurs 
and 
the first coordinate of $B_{i+1}$ is $b_{\ell}$ for some $\ell<n$
 is, with probability approaching $1$ as $n\rightarrow\infty$,
between $((1/18)-\epsilon^{\prime})n$ and $((1/18)+\epsilon^{\prime})n$.
This number of $i$ is within $1$ of the number of $a$ in 
Lemma~\ref{lemmathree}; the only possible difference occurs when 
$B_{i+1}$ has $b_n$ as one of its coordinates.
\qed

In conclusion, the number of $a$ in $\{1,\dots,n-1\}$ such that
$\tb_{a-1}=1$ and $\tb_a=1$ (given that we are in \firstone) lies,
with probability approaching $1$ as $n\rightarrow\infty$, between
$((4/18)-\epsilon)n$ and $((4/18)+\epsilon)n$ for each $\epsilon>0$.

\section{Stirling's Formula Argument}

If the number of $a$ in $\{1,\dots,n-1\}$ with $\tb_{a-1}=1$ and
$\tb_a=1$ is no more than $((4/18)+\epsilon)n$, then either the
number of odd $a$ in $\{1,\dots,n-1\}$ with $\tb_{a-1}=1$ and $\tb_a=1$
is no more than $((2/18)+\epsilon/2)n$
or the number of even $a$ in $\{1,\dots,n-1\}$ with $\tb_{a-1}=1$ and
$\tb_a$ is no more than $((2/18)+\epsilon/2)n$.

Let us suppose that $n$ is even and the number of odd $a$ in
$\{1,\dots,n-1\}$ with $\tb_{a-1}=1$ and $\tb_a=1$ is no
more than $((2/18)+\epsilon/2)n$ where $\epsilon>0$ is such that $(2/18)+
\epsilon/2<1/8$. Then, the number of possible values of $\sum_{i=0}^{n-1}
2^{n-1-i}b_i$ if we have \firstone is at most
\[
\sum_{(\ell_1,\ell_2,\ell_3,\ell_4)\in R_n}{((1/2)n)!\over\ell_1!\ell_2!
\ell_3!\ell_4!}
\]
where
$R_n=\{(\ell_1,\ell_2,\ell_3,\ell_4):\ell_1+\ell_2+\ell_3+\ell_4=(1/2)n,
\ell_1\le((2/18)+\epsilon/2)n\}$. The values $\ell_1$, $\ell_2$, $\ell_3$,
and $\ell_4$ represent the number of odd $a$ in $\{1,\dots,n-1\}$ such that
$\tb_{a-1}=1$ and $\tb_a=1$, $\tb_{a-1}=1$ and $\tb_a=0$,
$\tb_{a-1}=0$ and $\tb_a=1$, and $\tb_{a-1}=0$ and $\tb_a=0$, respectively. 

For some polynomial $p_1(n)$ of $n$,
\begin{eqnarray*}
&&\sum_{(\ell_1,\ell_2,\ell_3,\ell_4)\in R_n}{((1/2)n)!\over
\ell_1!\ell_2!\ell_3!\ell_4!}
\\
&\le&\left(\left({2\over 18}+{\epsilon\over 2}\right)n+1\right)n^2
{\left({1\over 2}n\right)!\over
\lfloor\left({2\over 18}+{\epsilon\over2}\right)n\rfloor!
\left(\lfloor\left({7\over 54}-{\epsilon\over 6}\right)n\rfloor!\right)^3}
\\
&\le&p_1(n){\left({1\over 2}n\right)^{(1/2)n}\over
\left(\left({2\over 18}+{\epsilon\over 2}\right)n\right)^{((2/18)+\epsilon/2)n}
\left(\left({7\over 54}-{\epsilon\over 6}\right)n\right)^{((7/54)-\epsilon/6)3n}
}\\
&=&p_1(n)2^{(0.5\log_2(0.5)-((2/18)+\epsilon/2)\log_2((2/18)+\epsilon/2)-
((7/18)-\epsilon/2)\log_2((7/54)-\epsilon/6))n}
\end{eqnarray*}

But if 
\[
c_1<{1\over 0.5\log_2(0.5)-{2\over 18}\log_2\left({2\over 18}\right)
-{7\over 18}\log_2\left({7\over 54}\right)}
\]
where $c_1$ is constant and $n=n(p)<c_1\log_2p$, 
then we can choose $\epsilon>0$ so that
\[
{p_1(n)2^{(0.5\log_2(0.5)-((2/18)+\epsilon/2)\log_2((2/18)+\epsilon/2)-
((7/18)-\epsilon/2)\log_2((7/54)-\epsilon/6))n}
\over p}\rightarrow 0
\]
as $p\rightarrow\infty$.

Thus for such values $n$, at most $o(p)$ values of
$\sum_{i=0}^{n-1}2^{n-1-i}b_i$ occur in \firstone if $n$ is even and the
number of odd $a$ in $\{1,\dots,n-1\}$ with $\tb_{a-1}=1$ and $\tb_a=1$ is no
more than $((2/18)+\epsilon/2)n$. Minor adaptations of this argument
apply if $n$ is odd, we consider the number of even $a$ instead of the number of
odd $a$, or we consider \firstnegone instead of \firstone. (For example,
if $n$ is odd
but we still consider odd $a$ in \firstone, note
that there are at most $o(p)$ different values of $X_{n-1}$ and $3$ different
values of $b_n$ to get that there are at most $o(p)$ different values of
$X_n$ in this case.)
There is only one value which is neither in \firstone nor in \firstnegone.

Observe that
\[
{1\over 0.5\log_2(0.5)-{2\over 18}\log_2\left({2\over 18}\right)
-{7\over 18}\log_2\left({7\over 54}\right)}
\approx 1.001525.
\]
Thus we may choose a value $c_1>1$ where if $n=n(p)<c_1\log_2p$, $X_n$ has, 
except for events with probability approaching $0$ as $p\rightarrow\infty$,
at most $o(p)$ values. Thus $\|P_n-U\|\rightarrow 1$ as $p\rightarrow\infty$.
\qed

\section{A Larger Value for $c_1$}

A more careful analysis of the proofs of Lemmas \ref{lemmaone}, 
\ref{lemmathree}, and \ref{lemmafour} 
shows that for each $\epsilon>0$,
the number of odd $a$
in $\{1,\dots,n-1\}$ with $\tb_{a-1}=1$ and $\tb_a=1$ and the number of
even $a$ in $\{1,\dots,n-1\}$ with $\tb_{a-1}=1$ and $\tb_a=1$ both lie
between $((2/18)-\epsilon)n$ and $((2/18)+\epsilon)n$ with probability
approaching $1$ as $n\rightarrow\infty$ given that we are in \firstone.
Since, given $j$, $C_a$ are independent when $a$ is odd, $a-1>j$, and
$a<n$ and $C_a$ are independent when $a$ is even, $a-1>j$, and $a<n$, 
the extension of Lemma~\ref{lemmafour} is straightforward.
To see how to extend Lemma~\ref{lemmaone}, consider the following argument.
Let $i_k$ be the $k$-th odd value of $i$ such that $D_i$ occurs, and let
$m_k$ be the value of $a$ such that $b_a$ is the first coordinate of the
block $B_{1+i_k}$. Note that
$m_2-m_1$, $m_3-m_2$, $m_4-m_3, \dots$ are i.i.d. Let $p$ be the
probability that $m_2-m_1$ is odd. If $m_1$ is even, then the sequence
$m_1, m_2, m_3, \dots$ consists of $r_1$ consecutive even values, then
$r_2$ consecutive odd values, then $r_3$ consecutive even values,
etc. where $r_1, r_2, r_3, \dots$ are i.i.d. geometric random
variables with parameter $p$. If $m_1$ is odd, then the sequence
$m_1, m_2, m_3, \dots$ consists of $r_1$ consecutive odd values,
then $r_2$ consecutive even values, then $r_3$ consecutive odd values,
etc. where $r_1, r_2, r_3, \dots$ are i.i.d. geometric random
variables with parameter $p$. Note that for each $\epsilon>0$, by 
Kolmogorov's maximal inequality (see p. 61 of Durrett~\cite{durrett},
for example), $\max(|r_1-r_2|, |(r_1-r_2)+(r_3-r_4)|, \dots,
|(r_1-r_2)+(r_3-r_4)+\dots+(r_{n-1}-r_n)|)<\epsilon n$ for even $n$ with
probability approaching $1$ as $n\rightarrow\infty$. Since for
some positive constant $c$, $\max(r_1, r_2, \dots, r_n)<c\ln(n)$ with 
probability approaching $1$ as $n\rightarrow\infty$, this result
and a similar result involving $D_i$ when $i$ is even
imply that, for each $\epsilon>0$, the number of odd $a$
in $\{1,\dots,n-1\}$ so that $b_{a-1}=1$, $b_a=1$, $\tb_{a-1}=1$,
and $\tb_a=1$ minus the number of even $a$ in $\{1,\dots,n-1\}$ with
$b_{a-1}=1$, $b_a=1$, $\tb_{a-1}=1$, and $\tb_a=1$ has absolute value less
than $\epsilon n$ with probability approaching $1$ as $n\rightarrow\infty$ 
given that we are in \firstone. Thus given $\epsilon>0$, the number of such
odd $a$ lies between $((1/36)-\epsilon)n$ and $((1/36)+\epsilon)n$ with
probability approaching $1$ as $n\rightarrow\infty$ given that we are in
\firstone. A similar argument applies for Lemma~\ref{lemmathree}.

With arguments resembling the proofs of Lemmas \ref{lemmaone}, 
\ref{lemmatwo}, \ref{lemmathree}, and \ref{lemmafour},
one can show

\begin{lemma}
\label{lemmafive}
Given that we are in \firstone, the number 
of $a$ in $\{1,\dots,n-1\}$
such that $\tb_{a-1}=1$ and $\tb_a=0$ lies, for each $\epsilon>0$, between
$((5/18)-\epsilon)n$ and $((5/18)+\epsilon)n$ with probability
approaching $1$ as $n\rightarrow\infty$.
\end{lemma}

\begin{lemma}
\label{lemmasix}
Given that we are in \firstone, the number 
of $a$ in $\{1,\dots,n-1\}$
such that $\tb_{a-1}=0$ and $\tb_a=0$ lies, for each $\epsilon>0$, between
$((4/18)-\epsilon)n$ and $((4/18)+\epsilon)n$ with probability
approaching $1$ as $n\rightarrow\infty$.
\end{lemma}

\begin{lemma} 
\label{lemmaseven}
Given that we are in \firstone, the number 
of $a$ in $\{1,\dots,n-1\}$
such that $\tb_{a-1}=0$ and $\tb_a=1$ lies, for each $\epsilon>0$, between
$((5/18)-\epsilon)n$ and $((5/18)+\epsilon)n$ with probability
approaching $1$ as $n\rightarrow\infty$.
\end{lemma}

While the details are not shown here, Table~\ref{thetable} outlines the
arguments to be shown. For example, the entry $1/18$ for $b_{a-1}=1$, 
$b_a=-1$, $\tb_{a-1}=0$, and $\tb_a=0$ means that the number of $a$ in
$\{1,\dots,n-1\}$ with $b_{a-1}=1$, $b_a=-1$, $\tb_{a-1}=0$, and $\tb_a=0$ 
lies, given $\epsilon>0$, between $((1/18)-\epsilon)n$ and
$((1/18)+\epsilon)n$ with probability approaching $1$ as $n\rightarrow\infty$
given that we are in \firstone. 

\begin{table}
\label{thetable}
\begin{tabular}{|l|c|c|c|c|}\hline
  &$\tb_{a-1}=0$, &$\tb_{a-1}=0$, & $\tb_{a-1}=1$, 
&$\tb_{a-1}=1$, \\
 &$\tb_a=0$&$\tb_a=1$&$\tb_a=1$&$\tb_a=0$ \\ \hline 
$b_{a-1}=1$, $b_a=1$&$0$ & $0$ & $1/18$ &$1/18$ \\ \hline
$b_{a-1}\ne 1$, $b_a\ne 1$&$1/9$ &$1/9$ &$1/9$&$1/9$\\ \hline
$b_{a-1}=0$, $b_a=1$& $1/18$& $1/18$&$0$&$0$\\ \hline
$b_{a-1}=-1$, $b_a=1$&$0$&$0$&$1/18$&$1/18$\\ \hline
$b_{a-1}=1$, $b_a=0$&$0$&$1/18$&$0$&$1/18$\\ \hline
$b_{a-1}=1$, $b_a=-1$&$1/18$&$1/18$&$0$&$0$\\ \hline
\end{tabular}
\caption{Cases for $\tb_{a-1}$, $\tb_a$, $b_{a-1}$, and $b_a$}
\end{table}

More careful arguments (similar to the extensions of Lemmas \ref{lemmaone},
\ref{lemmathree}, and \ref{lemmafour}) show that the number of odd $a$ in
$\{1,\dots,n-1\}$ with $\tb_{a-1}=1$ and $\tb_a=0$ lies between
$((5/36)-\epsilon)n$ and $((5/36)+\epsilon)n$ (given $\epsilon>0$)
with probability approaching $1$ as $n\rightarrow\infty$ given that we are
in \firstone. Similar statements hold for even $a$ here; similar
statements (where $2/18$ replaces $4/18$ and $5/36$ replaces $5/18$)
also hold for odd $a$ and even $a$ in Lemmas
\ref{lemmasix} and \ref{lemmaseven}. 

The total number of possible values of $\sum_{i=0}^{n-1}2^{n-1-i}b_i$
(except for events
with probability approaching $0$ as $n\rightarrow\infty$)
in \firstone is at most (for even $n$)
\[
\sum_{(\ell_1,\ell_2,\ell_3,\ell_4)\in S_n}{{1\over 2}n\choose 
\ell_1,\ell_2,\ell_3,\ell_4}
\]
with $S_n=\{(\ell_1,\ell_2,\ell_3,\ell_4):\ell_1+\ell_2+\ell_3+\ell_4=(1/2)n,
((4/36)-\epsilon)n<\ell_1<((4/36)+\epsilon)n,
((5/36)-\epsilon)n<\ell_2<((5/36)+\epsilon)n,
((5/36)-\epsilon)n<\ell_3<((5/36)+\epsilon)n,
((4/36)-\epsilon)n<\ell_4<((4/36)+\epsilon)n\}$.

A Stirling's formula argument shows that if
\[
c_1<{1\over 0.5\log_2(0.5)-{4\over 18}\log_2\left({4\over 36}\right)
-{5\over 18}\log_2\left({5\over 36}\right)}
\]
and $n=n(p)<c_1\log_2p$ where $c_1$ is a constant, then
\[
\sum_{(\ell_1,\ell_2,\ell_3,\ell_4)\in S_n}{{1\over 2}n\choose 
\ell_1,\ell_2,\ell_3,\ell_4}
\]
is $o(p)$. For odd $n$ or \firstnegone, similar arguments can
be used. Thus if $n=n(p)<c_1\log_2p$, $X_n$ has
$o(p)$ possible different values except for events with probability
approaching $0$ as $p\rightarrow\infty$. Thus $\|P_n-U\|\rightarrow 1$
as $p\rightarrow\infty$.

Note that
\[
{1\over 0.5\log_2(0.5)-{4\over 18}\log_2\left({4\over 36}\right)
-{5\over 18}\log_2\left({5\over 36}\right)}\approx 1.00448.
\]
Thus there is a gap between this 
lower bound and the best upper bound claimed
in Chung, Diaconis, and Graham~\cite{cdg}. Exploring this gap is a potential
problem for further study.

\section{Acknowledgments}

The author thanks Ron Graham for mentioning this problem
in a talk at a conference on the mathematics of Persi Diaconis in 2005.
The author also thanks Persi Diaconis for encouragement.

\end{document}